\documentclass[12pt]{article}
\usepackage{amsbsy,amsfonts,amsmath,amssymb,amsthm,bbm,enumerate,euscript,graphicx}
\newtheorem{Theorem}{Theorem}[section]
\newtheorem{Proposition}[Theorem]{Proposition}
\newtheorem{Lemma}[Theorem]{Lemma}

\newtheorem{Corollary}[Theorem]{Corollary}

\numberwithin{equation}{section}

\newcommand{\tPo}{{\cal H}}

\newcommand{\bpr}{\begin{proof}}
\newcommand{\epr}{\end{proof}}
\newcommand{\0}{{\bf o}}
\newcommand{\bfe}{{\bf e}}
\newcommand{\bx}{{\bf x}}
\newcommand{\cc}{{\rm c}}
\newcommand{\bdm}{\begin{displaymath}}
\newcommand{\edm}{\end{displaymath}}

\renewcommand{\H}{{\cal H}}
\newcommand{\cH}{{\cal H}}

\renewcommand{\emptyset}{\varnothing}

\newcommand{\bean}{\begin{eqnarray*}}
\newcommand{\eean}{\end{eqnarray*}}
\newcommand{\bea}{\begin{eqnarray}}
\newcommand{\eea}{\end{eqnarray}}
\newcommand{\be}{\begin{eqnarray}}
\newcommand{\ee}{\end{eqnarray}}
\newcommand{\benu}{\begin{enumerate}}
\newcommand{\eenu}{\end{enumerate}}

\newcommand{\toP}{\stackrel{P}{\longrightarrow}}

\newcommand{\rem}{\paragraph{Remark.}}

\def\E{{\mathbb{E}}}
\def\Pr{{\mathbb{P}}}

\def\R{\mathbb{R}}

\def\Z{\mathbb{Z}}

\def\1{{\bf 1}}
\def\N{\mathbb{N}}

\def\Po{{\cal P}}

\def\eps{{\varepsilon}}

\def\X{{\cal X}}

\def\dist{{\rm dist}}

	\newcommand{\cC}{{\cal C}}
	\newcommand{\cV}{{\cal V}}

 \newcommand{\la}{\lambda}

 {\nopagebreak {\hspace*{\fill}\rule{7pt}{7pt}}\\ }

\title{Giant component of the soft random geometric graph
}

\author{Mathew D. Penrose$^1$\\
University of Bath
}
 \footnotetext{ Department of
Mathematical Sciences, University of Bath, Bath BA2 7AY, \\
United
Kingdom: {\texttt m.d.penrose@bath.ac.uk} }

\begin{document}
\maketitle

\begin{abstract}
	Consider a 2-dimensional soft random geometric graph $G(\lambda,s,\phi)$, obtained by placing a Poisson($\lambda s^2$) number of vertices uniformly at random in a square of side $s$, with  edges placed  between each pair $x,y$ of vertices with probability $\phi(\|x-y\|)$, where $\phi: \R_+ \to [0,1]$ is a finite-range connection function.  This paper is concerned with the asymptotic behaviour of the graph $G(\lambda,s,\phi)$ in the large-$s$ limit with $(\lambda,\phi)$ fixed.  We  prove that the proportion of vertices in the largest component converges in probability to the  percolation probability for the corresponding random connection model, which is a random graph defined similarly for a Poisson process on the whole plane. We do not cover the case where $\lambda$ equals the critical value $\lambda_c(\phi)$.

AMS classifications:
60C05, 60D05, 60K35.
\end{abstract}

\section{Introduction and statement of results}

%\label{ch020723b}

%\subsection{The random connection model}
%\label{secconper}
%\index{percolation!continuum}
Let $\phi: [0,\infty) \to [0,1]$ be a 
nonincreasing function.
Given a locally finite point set $\X \subset \R^2$, 
%p\in (0,1]$,
let
$G(\X,\phi)$ denote the graph with vertex set $\X$, where
for each $\{x,y\} \subset \X$, the
edge $xy$ is included with probability 
%$p{\bf 1}_{\{\||x-y\| \leq 1\}}$,
$\phi (\|y-x\|)$,
independently of
the other pairs. Here $\|\cdot\|$ is the Euclidean norm.

Given $\la >0$, 
let $\H_\la$ denote a homogeneous Poisson process
of intensity $\la$ in $\R^2$.
Given also $s >0$, set $B(s):=[-s/2,s/2]^2$ and let $\H_{\la,s}$ denote the restriction
of $\H_\la$ to $B(s)$, which is
 a homogeneous Poisson process of intensity $\la$ in $B(s)$.
We are interested in the graphs $G(\H_\la,\phi)$ and 
$G(\H_{\la,s},\phi)$, which are
known as the  {\em random connection model} \cite{MR}
and {\em soft random geometric graph} \cite{P2}
respectively, with connection function $\phi$.

Let $\tPo^{\0}_{\lambda}$ denote the point process $\tPo_{\lambda} \cup \{\0\}$,
where $\0$ is the origin in $\R^2$.
For $k \in \N$, let $\pi_k(\phi,\lambda)$ denote the probability that
the component of $G(\tPo^{\0}_{\lambda},\phi)$   containing  the origin is of
order $k$.
%see (\ref{per2a}) below for a  formula for $p_k(\lambda)$.
 The {\em percolation probability}
\index{percolation probability} $\theta(\phi,\lambda)$ is
 the probability that $\0$ lies in an
infinite component of the graph $G(\tPo^{\0}_{\lambda},\phi)$, that
is, 
%and is defined  by 
$$
\theta (\phi,\lambda) : = 1 - \sum_{k=1}^ \infty \pi_k(\phi,\lambda).
$$
A standard coupling argument shows that $\theta(\phi,\la)$ is nondecreasing in $\la$.
%\begin{Exercise}
	%{\rm
	%Suppose $0< \lambda < \lambda'$. Show that
	%$p_\infty(\lambda) \leq p_\infty(\lambda')$.}
%\end{Exercise}
The {\em critical value}\index{critical value}
 (continuum percolation threshold)
\index{continuum percolation threshold}
  $\lambda_\cc(\phi)$  % that is,
  is defined by
\bea
\lambda_\cc(\phi) := \inf \{\lambda>0: \theta(\phi,\lambda)>0\}.  
\label{lambdacdef}
\eea
It is known (see \cite{MR}) that 
	$0 < \lambda_\cc < \infty$, provided $0 < \int_{\R^2} \phi(\|x\|)dx
	< \infty$.

	\newpage

For any finite graph $G$, let $L_j(G)$  denote the order of
its $j$th-largest component, that is, the $j$th-largest of the 
orders of its components, or zero if it has fewer than $j$
components. In this paper we prove the following results about 
convergence in probability of $s^{-2}L_j(G(\H_{\la,s},\phi))$ for
$j=1,2$, as $s \to \infty$ with $\phi ,\la$ fixed.
%in Theorems \ref{thgiantsub} and \ref{thm1}.

\begin{Theorem}
	\label{thgiantsub}
	Suppose %$p \in (0,1]$ and
	$\lambda >0$ with $\theta(\phi,\la)=0$.
   Then
\bea
s^{-2} L_1 ( G(\tPo_{\lambda,s},\phi) ) \toP  0
{\rm ~~~~~~as~~} s \to \infty.
\label{eqsub}
\eea
\end{Theorem}

\begin{Theorem}
	\label{thm1}
	Suppose $\sup\{r >0: \phi(r) >0\} \in (0,\infty)$
	%$p \in (0,1] $
	and $\la > \la_c(\phi)$. Then
	as $s \to \infty$ we have that
	$s^{-2}L_{1}(G(\H_{\la,s},\phi)) \toP \la \theta(\phi,\la)$ and 
	$s^{-2}L_{2}(G(\H_{\la,s},\phi)) \toP 0.$
\end{Theorem}

These theorems  do not address the case with $\lambda= \lambda_c(\phi)$,
unless we know $\theta(\phi,\la_c(\phi)) =0$ (if we had this, we could
apply Theorem \ref{thgiantsub} when $\lambda = \lambda_c(\phi)$). 
Theorem \ref{thm1} also does not address the case of $\phi$ having unbounded
range (recall we are assuming $\phi$ is nonincreasing).

In the case with $\phi = {\bf 1}_{[0,1]}$, these results
were already proved in \cite{Penbk}, but the method here provides
an alternative and possibly shorter proof (the proof in \cite{Penbk}
relies on a lengthy RSW argument from \cite{MR}, as well as not
working for general $\phi$).
When $\phi= {\bf 1}_{[0,1]}$ it is known \cite{MR}
that $\theta(\phi,\la_c(\phi))=0$ since this case is equivalent to a Boolean
model.

%	\begin{Lemma}
%		As a function of $\lambda$, the percolation
%	probability  $\theta(\lambda)$ is right continuous. Also 
%		$\theta(\lambda)$ is {\em continuous}
%		in $\la$, except possibly at $\la = \la_c(\phi)$.
%	\end{Lemma}
%Moreover,
%it is known \cite{LPZ} that $p_\infty(\lambda)$ is infinitely differentiable
% in $\lambda$, except at $\lambda = \lambda_c$.

\section{Proof of theorems}
\label{chgiant1}
\label{seccoper}

\subsection{Preliminaries}

We introduce some further notation that will be used in proving
Theorems \ref{thgiantsub} and \ref{thm1}.

For $\X \subset \R^2$ and $x,y \in \R^2$ we sometimes write
$\X^x$ for $\X \cup \{x\}$ and $\X^{x,y}$ for $\X \cup \{x,y\}$.

Let $\phi: [0,\infty) \to [0,1]$ be nonincreasing.
%\in (0,1]$.
We shall view $\phi$ as fixed from now on, and
for any locally finite $\X \subset \R^2$ 
write simply $G(\X)$ instead of $G(\X,\phi)$.
Also for
 $x \in \X$,
let $\cC_x(\X)$ denote the vertex set of the
component of $G(\X)$ containing $x$.

Given $x \in \R^2$ and $r >0$,
we write $D_r(x)$ for the disk $\{y \in \R^2:\|y-x\| \leq r\}$
and $D_r$ for $D_r(\0)$.  Also let $S_r:= B(2r) = [-r,r]^2$, and let $\bfe := (1,0)$. 

Given $A,B \subset \R^2$, and locally finite $\X \subset \R^2$,
we write $\{ A \leftrightarrow B $ in $G(\X)\}$
for the event that there exist $x \in A \cap \X$ and $y \in B \cap \X$
such that there is a path in $G(\X)$ from $x $ to $y$.
We write $\{ A \leftrightarrow \infty $ in $G(\X)\}$
for the intersection over all $n \geq 1$ of events
 $\{ A \leftrightarrow \R^2 \setminus D_n $ in $G(\X)\}$

Next, we assemble some known facts which will be used later.

We say a real-valued function $f$, defined on graphs
$(\cV,E)$,  with
$\cV \subset \R^2$ locally finite, is {\em increasing}  if $f(\cV,E) \leq
f(\cV',F)$
whenever $\cV \subset \cV'$ and $E \subset E'$. We say $f$ is {\em decreasing}
if $-f$ is increasing. Given $\la >0$ and given $\phi$,
we say $E$ is an increasing (resp. decreasing) event
on $G(\H_\la)$ if ${\bf 1}_E $ is an increasing (resp. decreasing)
function of $G(\H_\la)$.
\begin{Lemma} [Harris-FKG inequality]
	\label{lemFKG}
	Suppose $f,g$ are   measurable bounded  increasing real-valued
	functions
	defined on graphs $(\cV,E)$ with $\cV \subset \R^2$ locally finite,
%	point configurations in $\R^2$. 
	%{\bf (glossing over what measurable means!)}. 
	%Let $\la >0$.
	Then $$\E[f(G(\H_{\la})) g(G(\H_{\la})) ] \geq \E[f(g(\H_{\la}))]
	\E[g(G(\H_{\la}))]. $$
	The same inequality holds if $f$ and $g$ are both decreasing.
\end{Lemma}
\bpr
See \cite{HLM}, where measurability issues are also dealt with.
\epr

\begin{Corollary}[Square Root trick]
	Let $\lambda >0$, $k \in \N$, $\eps \in (0,1)$. 
	Suppose for $i=1,\ldots,k$ we
	have  increasing events $A_i$
	defined on $G(\H_{\la})$,
	such that $\Pr[ \cup_{i=1}^k A_{i}] > 1 - \eps$.

	Then $ \max_{1 \leq i \leq k} \Pr[A_{i}] > 1- \eps^{1/k}$.
\end{Corollary}
\bpr
Set $M = \max_{1 \leq i \leq k} \Pr[A_{i}]$. The
events $A_{i}^c$ are all decreasing, so by Lemma \ref{lemFKG},
%the Harris-FKG
%inequality
\bean
\eps >
\Pr[ \cap_{i=1}^k A_i^c ] \geq \prod_{i=1}^k \Pr[A_i^c]
\geq (1-M)^k,
\eean
so that $1-M < \eps^{1/k}$ and $M > 1- \eps^{1/k}$.
\epr

%\section{Uniqueness of the infinite component}
%Given $p \in (0,1), \la >0 $,
Given $\la >0 $,
let $N_\infty(\phi,\la)$ be the number of infinite components of the graph
$G(\H_\la) $. It is not hard to show that
	if $\theta(\phi,\lambda)=0$,
	then $\Pr[N_\infty(\phi,\la) =0] =1$.
The next preliminary result concerns {\em uniqueness of the infinite cluster},
in
 the other case,
where $\theta(\phi,\lambda) >0$.
%we have a deeper result.
%`uniqueness of the infinite cluster':

\begin{Lemma}
\label{perthuni}
	Suppose $\theta(\phi,\lambda) > 0$. Then $\Pr[N_\infty(\phi,\la)=1]=1$.
\end{Lemma}
\bpr
See \cite{MR}.
\epr
Another useful fact is the {\em Mecke formula} for the random
connection model. Let $s, \la >0$ and suppose $f(x,G) \in \R_+$ is defined
for
all pairs $(x,G)$ where $G$ is a finite graph with vertex
set $\cV(G) \subset B(s)$ and $x \in \cV(G)$. Then whenever
the following expectations are defined we have
\bea
\E \sum_{x \in \H_{\la,s}} f(x,G(\cH_{\la,s})) =
\la \int_{B(s)} \E[f(x, G( \cH_{\la,s}^x))] dx
\label{eqMecke1}
\eea
and moreover if $g(x,y,G) \in \R_+$ is defined whenever additionally
$y \in \cV(G)$ then
\bea
\E\sum \sum_{x,y \in \H_{\la,s}, x \neq y}
g(x,y,G(\cH_{\la,s})) = \la^2 \int_{B(s) } \int_{B(s)}
\E[ g(x,y,G(\cH_{\la,s}^{x,y})) ] dy dx.
\label{eqMecke2}
\eea
The Mecke formulae (\ref{eqMecke1}) and (\ref{eqMecke2})
can be derived by conditioning on $\cH_{\la,s}$ and
using the usual Mecke formulae from e.g. \cite{LP}.

Also of use to us is the following {\em  sequential construction}
of clusters in $G(\cH_{\la,s})$. Let $\la, s>0$ and $A \subset B(s)$
(typically a disk).
The set $\cup_{x \in \H_{\la} \cap A} \cC_x(\cH_{\la,s})$
can be created as follows:

First generate $\H_\la \cap A$. Denote the points so created 
as {\em active points} and let the initial  intensity function
of {\em unexplored  points}, i.e. Poisson points that
are not yet generated, is $g = \la {\bf 1}_{B(s) \setminus A} (\cdot)$.

Next, choose an active point $x$ and generate a Poisson process
of intensity 
%$p g(\cdot) {\bf 1}_{D_1(x)} (\cdot)$,
$ g(\cdot) \phi(\cdot -x) $,
representing the
previously unexplored points
of $\cH_{\la,s}$ 
that are connected directly to $x$. Label all the new 
points as `active', and change the status
of $x$ from `active' to `finished'. Also change the intensity of
unexplored points from $g(\cdot)$ to
$g(\cdot) (1- \phi(\cdot -x))$.
%$g(\cdot) [ {\bf 1}_{\R^2 \setminus D_1(x)}(\cdot)
%	+ (1-p) {\bf 1}_{D_1(x)}(\cdot) ]$.

Then pick a new active point and repeat the above, using the new intensity of
unexplored points. Keep repeating until we run out of active points,
then stop.

We shall refer to the above procedure as
{\em growing the cluster sequentially}. This method is described in detail
for the case $\phi = {\bf 1}_{[0,1]}$ %$p=1$
in \cite{P96}, and for the general random connection
model in \cite{MPS}.

%We can  do the same thing using non-constant Poisson intensities,
%for example if $A \subset B(s)$
%we can generate $\cup_{x \in A} \cC_x(\cH_{\la,s})$ the same way
%but with the initial intensity
%of unexplored points set to $g(\cdot) = {\bf 1}_{B(a) \setminus A}(\cdot)$.

%\allco

\subsection{The subcritical case}

\bpr[Proof of Theorem \ref{thgiantsub}]
Suppose $\la > 0$ with  $\theta(\phi,\lambda)=0$.

Let $\eps >0$. Let $N_s := \sum_{x \in \H_{\la,s}} {\bf 1}
\{|\cC_x(\H_{\la,s})| \geq \eps s^2\}$.
%be the number of $x \in \H_{\la,s}$ such that
%$|\cC_x(\H_{\la,s})| \geq \eps s^2$.
If $L_1(G(\H_{\la,s})) \geq \eps s^2$,
then $N_s  \geq \eps s^2$. Hence by Markov's inequality and the Mecke formula,
\bean
\Pr[ L_1(G(\H_{\la,s})) \geq \eps s^2 ]
\leq (\eps s^2)^{-1} \E[N_s]
= (\eps s^2)^{-1} \int_{B(s)} \Pr[|\cC_x(\H_{\la,s}^x ) | \geq \eps s^2] \la dx
\\
\leq (\eps s^2)^{-1} \int_{B(s)} \Pr[|\cC_x(\H_{\la}^x ) | \geq \eps s^2] \la dx
\\
= \la \eps^{-1} \sum_{k \geq \eps s^2} \pi_k(\phi,\lambda) 
\eean
which tends to zero as $s \to \infty$. Therefore
$s^{-2} L_1(G(\H_{\la,s})) \toP 0 $.
%converges in probability to zero.
%and hence so does $s^{-2} L_2(G(\H_{\la,s},1)) $. 
\epr

 \subsection{Renormalization}

 From now on we assume $\phi$ is nonincreasing with
 $\inf\{r>0: \phi(x) =0\} =1$. We shall prove Theorem \ref{thm1}
 only for this case, since simple scaling arguments then yield the general
 finite-range case in the statement of the theorem.

 Given  $\la, K, L, M \in (0,\infty)$ with $L >K, M> 2K$, define the following events:
 \begin{itemize}
		 \item
			 $U_{K,L,\la}$ is the event that
			 there is a unique component
			 of $G(\H_\la \cap D_{L+1})$ that meets
			 both $D_K$ and $\R^2 \setminus D_L$.
		 \item
			 $F_{K,M,\la} =\{ 
			 D_K  \leftrightarrow  D_K(M\bfe)$ in
			 $G(\H_\la \cap D_{3M}) \} $.
 \end{itemize}
 \begin{Proposition}
	 \label{renormprop}
	 Suppose $\la > \la_c(\phi) $ and
	 let $\eps \in (0,1)$.
	 There exist finite constants $K>0$ and $M > 3K$ such
	 that 
	 (i) $\Pr[D_K \leftrightarrow \infty$ in $G(\cH_{\la}) ] > 1-\eps$,
	 and (ii)
	 $\Pr[U_{K,M/3,\la} ]> 1-\eps $, and (iii)
	 $\Pr[F_{K,M,\la} ] > 1-\eps$.
 \end{Proposition}
 We shall use this to establish the
 limiting behaviour of $s^{-2}L_1(G(\H_{\la,s}))$ for
 $\la > \la_c(\phi)$.
 The point is that we can use it to compare $G(\H_\la)$
 with a finite-range dependent percolation process on the lattice $ \Z^2$.

 To prepare for the proof of Proposition \ref{renormprop}, fix $\la > \la_c(\phi)$ and $\eps \in (0,1)$.
 Set $\mu= (\la_c(\phi) + \la)/2$, so
 $\mu \in (\la_c(\phi),\la)$.
It is not hard to see that we can (and do) choose $\eta >0$ such that for
any two distinct points $x,y \in [0,3]^2$,
\bea
\Pr[ \{x \} \leftrightarrow \{y \} ~{\rm in}~ 
G(\H_{\la - \mu}^{x,y} \cap [0,3]^2) ]\geq \eta.
\label{etadef}
\eea
We then choose
$\nu \in \N$ such that $(1-\eta)^{\nu/9} < \eps/3$.
Then take
  \bea
  \delta : = (\eps/3) e^{-25 \la \nu}.
  \label{deltadef}
  \eea

 \bpr[Proof of Proposition \ref{renormprop}]
 We adapt  an argument in \cite{DCST}.
 Let $\mu, \eta, \nu$ and $\delta $  be as given above and
 let $\eps_1  = (\delta/3)^{32}$.
 Choose $K$ such that
 $\Pr[D_K \leftrightarrow \infty {\rm ~in~} G(\H_\mu) ] > 1- \eps_1$.
 Since $\eps_1< \eps$ and $\mu < \la$ this yields (i) at once.
 %Since $\delta > \eps_1$,
 Since $\eps_1 < \eps$, we
 claim that we can (and do) choose $n_1 \in \N$ with $n_1 > K$ such
 that  
 \bea
 \Pr[U_{K,n,\la}] > 1-  \eps, ~~~ \forall n \in [n_1,\infty).
 \label{n1eq1}
 \eea
 We leave the proof of this claim, 
	  using Lemma \ref{perthuni}, as an exercise.

 Now for integer $n \geq n_1$,   and
 for $0 \leq \alpha \leq \beta \leq n+1$, define the  event
 $$
 E_n (\alpha,\beta):= \{ 
 D_K \leftrightarrow [n,n+1] \times [\alpha, \beta] 
 ~{\rm in}~ G(\H_\mu \cap S_{n+1}) \}.
 $$
 Since $\Pr[D_K \leftrightarrow  S_n^c ~{\rm in}~G(\H_\mu \cap S_{n+1}) ] \geq
 \Pr[D_K \leftrightarrow \infty~{\rm in} ~G(\H_\mu)] > 1-\eps_1$,
 using the Square Root trick we can deduce that
 \bea
 \Pr[E_n(0,n+1)] > 1 - \eps_1^{1/8}.
 \label{0309a}
 \eea
 Next, observe that for fixed $n$, we have that
 as a function of $\alpha$,
 $\Pr[E_n(0,\alpha) ] -
 \Pr[E_n(\alpha,n+1)] $
 increases continuously from a value of
  $- \Pr[E_n(0,n+1)]$ at $\alpha =0$ to
  a value of
  $+ \Pr[E_n(0,n+1)]$ at $\alpha =n+1$.
 Therefore we can and do choose
 $\alpha_n \in (0,n+1)$ such that 
 $
 \Pr[E_n(0,\alpha_n) ]= \Pr[E_n(\alpha_n,n+1)].
 $
 %\begin{Exercise} 
%	 Prove the various assertions in the preceding lines.
% \end{Exercise} 
 Since $E_n(0,n+1) = E_n(0,\alpha_n) \cup E_n(\alpha_n,n+1)$,
 by (\ref{0309a}) and a further application of the Square Root trick
 we obtain that 
 \bea
 \Pr[E_n(\alpha_n,n+1)] =
 \Pr[E_n(0,\alpha_n)] > 1- \eps_1^{1/16}.
 \label{DCTS2}
 \eea

 By yet another application of the Square Root trick we obtain that
 $$
 \max (
 \Pr[E_n(0,\alpha_n/2)], 
 \Pr[E_n(\alpha_n/2,\alpha_n) ] ) > 1- \eps_1^{1/32} 
 $$
 so we can and do choose $y_n$,
  with either $y_n = \alpha_n/4$ or $y_n = 3 \alpha_n/4$,
 such that 
 \bea
 \Pr[E_n(y_n - \alpha_n/4,y_n + \alpha_n/4)] > 1 - \eps_1^{1/32}.
 \label{DCTS3}
 \eea

 Set $n_2 = 3n_1$.
 We claim that there exists integer $N \geq  n_2$ such that
 $\alpha_{3N} < 4 \alpha_N$.
 Indeed, if this were not true then  we would have
 for all $k \geq 1$ that
 $\alpha_{3^k n_2} \geq 4^k \alpha_{n_2}$, but since $\alpha_n \leq n+1$
 for all $n$, this would imply
 $3^k( n_2 +1) \geq 4^k \alpha_{n_2}$ so that $(4/3)^k \leq 
 (n_2+1)/\alpha_{n_2}$
 for all $k$, which is not true (since $\alpha_{n_2} >0$),
 justifying the claim.

 Choose (deterministic) integer $N \geq n_2 $ such that 
 $\alpha_{3N } < 4 \alpha_N$. 
Then
 by (\ref{DCTS2}) and (\ref{DCTS3}), setting $\eps_2:= \eps_1^{1/32}$
 we have
 \bea
 \min ( \Pr[E_N(\alpha_N,N+1)],
 \Pr[E_{3N}(y_{3N} - \alpha_{3N}/4, y_{3N} + \alpha_{3N}/4)] ) > 
 1- \eps_2.
 \label{0325a}
 \eea

 Now set $\bx = (2N,y_{3N})$ (we use bold face to indicate certain fixed 
 2-vectors such as $\0$ and $\bfe$).
 Let
  $S_{N+1}(\bx) := S_{N+1} + \bx = [N-1,3N+1] \times [y_{3N}-N -1,
  y_{3N} + N +1 ]$.
 Define the vertical blocks (see Figure 1)
 $$
 I := [3N,3N+1] \times [y_{3N}-\alpha_{3N}/4, y_{3N} + \alpha_{3N}/4],
 $$
  $$J^+ := [3N,3N+1] \times [y_{3N}+\alpha_N,y_{3N}+N+1],
   $$
   $$
   J^- := [3N,3N +1] \times [y_{3N}-N-1, y_{3N}- \alpha_N].
   $$

  Let $A^+$ be the event 
  $\{ D_{K}(\bx) \leftrightarrow J^+$ in
  $G(\H_\mu \cap S_{N+1}(\bx))\}$,
  and
  let $A^-$ be the event 
  $\{D_{K}(\bx) \leftrightarrow J^-$ in 
  $G(\H_\mu \cap S_{N+1}(\bx))\}$.
  Then $\Pr[A^+] = \Pr[A^-]= \Pr[E_N(\alpha_N,N+1)]$.
 Also $ E_{3N}(y_{3N}-\alpha_{3N}/4,y_{3N}+ \alpha_{3N}/4) =
  \{D_K \leftrightarrow I$ in $G(\cH_\mu \cap S_{3N+1})\}$.
  By %(\ref{n1eq1}),  
  (\ref{0325a}) and the union bound, 
  $$
  \Pr[ %U'_{K,n_1,\la} \cap
  A^+ \cap A^- \cap
  %E_{3N}(y_{3N}-\alpha_{3N}/4,y_{3N}+ \alpha_{3N}/4)
  \{D_K \leftrightarrow I {\rm ~ in ~} G(\cH_\mu \cap S_{3N+1})\}
  ]
  > 1-3\eps_2 = 1- \delta.
  $$

 Set $M = \|\bx\|$. Then $M \geq 2N \geq 3n_1$. By (\ref{n1eq1}),
 $ \Pr[U_{K,M/3,\la}] > 1- \eps$ so we have (ii).
  The proof is then completed by
  the following `gluing lemma'. \epr
  \begin{Lemma}
	  \label{lemglue}
	  If $\Pr[
  A^+ \cap A^- \cap 
	  %E_{3N}(y_{3N}-\alpha_{3N}/4,y_{3N}+ \alpha_{3N}/4)]
  \{D_K \leftrightarrow I {\rm ~ in ~} G(\cH_\mu \cap S_{3N+1})\}
  >1- \delta$, then $\Pr[F_{K,M,\la}] > 1- \eps$, where we take $M= \|\bx\|$.
  \end{Lemma}
  \rem{The proof below is not needed  for the special case
  $\phi ={\bf 1}_{[0,1]}$, since in this case
 % when $p=1$,
  the lemma is immediate because 
  $A^+ \cap A^- \cap 
  \{D_K \leftrightarrow I {\rm ~ in ~} G(\cH_\mu \cap S_{3N+1})\}$
	  implies $F_{K,M,\mu}$. 
	  }
\bpr[Proof of Lemma \ref{lemglue}]
  Divide $\R^2$ into half-open rectilinear
  squares $Q_i$ of side 1 and for each $i$ let $Q_i^+$ (respectively
  $Q_i^{++}$
  be the half-open square of side 3 (resp. 5)  with the same centre.
  We shall define a random variable $Z$ taking values in
  $\Z_+ \cup \{+\infty\}$, as follows.
  
  Grow the cluster $\cC := \cup_{x \in \cH_\mu \cap D_K} 
  \cC_x(\cH_\mu \cap S_{3N+1})$ sequentially.
  Let $\Po$ be the point process of unexplored points of 
  $\cH_\mu \cap S_{N+1}(\bx)$ at
  the end of this procedure, i.e.
  the points of $\H_\mu \cap S_{N+1}(\bx)$ that
  either lie outside $S_{3N+1}$ (since $S_{N+1}(\bx)$
  is conceivably not entirely contained in $S_{3N+1}$) or are 
  not connected by an edge to any point of $\cC$.
  If $\cC \cap D_K(\bx) \neq \emptyset$ then set $Z= +\infty$.

  Next, assuming $Z < \infty$, 
  grow the cluster
  $\cup_{x \in \Po \cap D_K(\bx)} \cC_x(\Po)$
  sequentially
  but do not continue the exploration from any points created
  that lie in $\cup_{\{i:\cC \cap Q_i \neq \emptyset\}}Q_i^+$; leave
  these points as `active'.
  We denote this second cluster by $\cC'$.
  Let ${\cal I} := \{i: \cC \cap Q_i \neq \emptyset, \cC' \cap Q_i^+
  \neq \emptyset\}$ and let $Z= |{\cal I}|$, as illustrated in Figure 1.

  Define the event 
  $E := \{D_K \leftrightarrow I {\rm ~ in ~} G(\cH_\mu \cap S_{3N+1})\} =
  \{\cC \cap I \neq \emptyset\}$.

  %Write just $E$ for
  %$\{D_K \leftrightarrow I {\rm ~ in ~} G(\cH_\mu \cap S_{3N+1})\}$
  %(this is the same as the event $\{\cC \cap I \neq \emptyset\}$).
  Recall that $\nu$ was defined just after (\ref{etadef}).
Suppose $\nu < Z < \infty$.
Then we can find $i_1, \ldots i_{\lceil \nu/9 \rceil}
\in {\cal I}$ such that the squares
$Q_{i_1}^+, \ldots, Q^+_{i_{\lceil \nu/9 \rceil}}$
are disjoint.
Then sprinkling an independent  Poisson process $\H_{\la -\mu}$
on top of $\H_\mu$ in these squares, for each square we have a chance
at least $\eta$ to join $\cC$ to $\cC'$ via the sprinkled points in
that square. Since we may assume $\cH_\la=
\cH_\mu \cup \cH_{\la - \mu}$, 
setting $F'= \{D_K \leftrightarrow D_K(\bx) ~{\rm in} ~ G(\H_\la
\cap (S_{3N+1} \cup S_{N+1}(\bx))) \}$
we obtain that
\bea
\Pr[  (F')^c
| E \cap \{Z > \nu\}] \leq (1- \eta)^{\nu/9} < \eps/3.
\label{0406a}
\eea
Note that if $Z=+\infty$ then $F'$ must occur so this case
is included in (\ref{0406a}).

  Now suppose $Z \leq \nu$. If also event $E$ occurs
  then it is not possible to find  paths in $S_{N+1}(\bx)$ both
  from $D_K(\bx)$ to $J^+$, and from $D_K(\bx) $ to $J^-$,
  with neither path passing through $\cup_{\{i: \cC \cap Q_i \neq \emptyset\}}
  Q_i^+$
  (see Figure 1, and also \cite[Figure 2.2]{DCST}; note
   $I \cap J^+ = I \cap J^- = \emptyset$ since
  $\alpha_{3N}/4 < \alpha_N$).
  Therefore we have not yet achieved event $A^+ \cap A^-$
  at this stage.

  \begin{figure}[!h]
\label{fig0}
\center
\includegraphics[width=10cm]{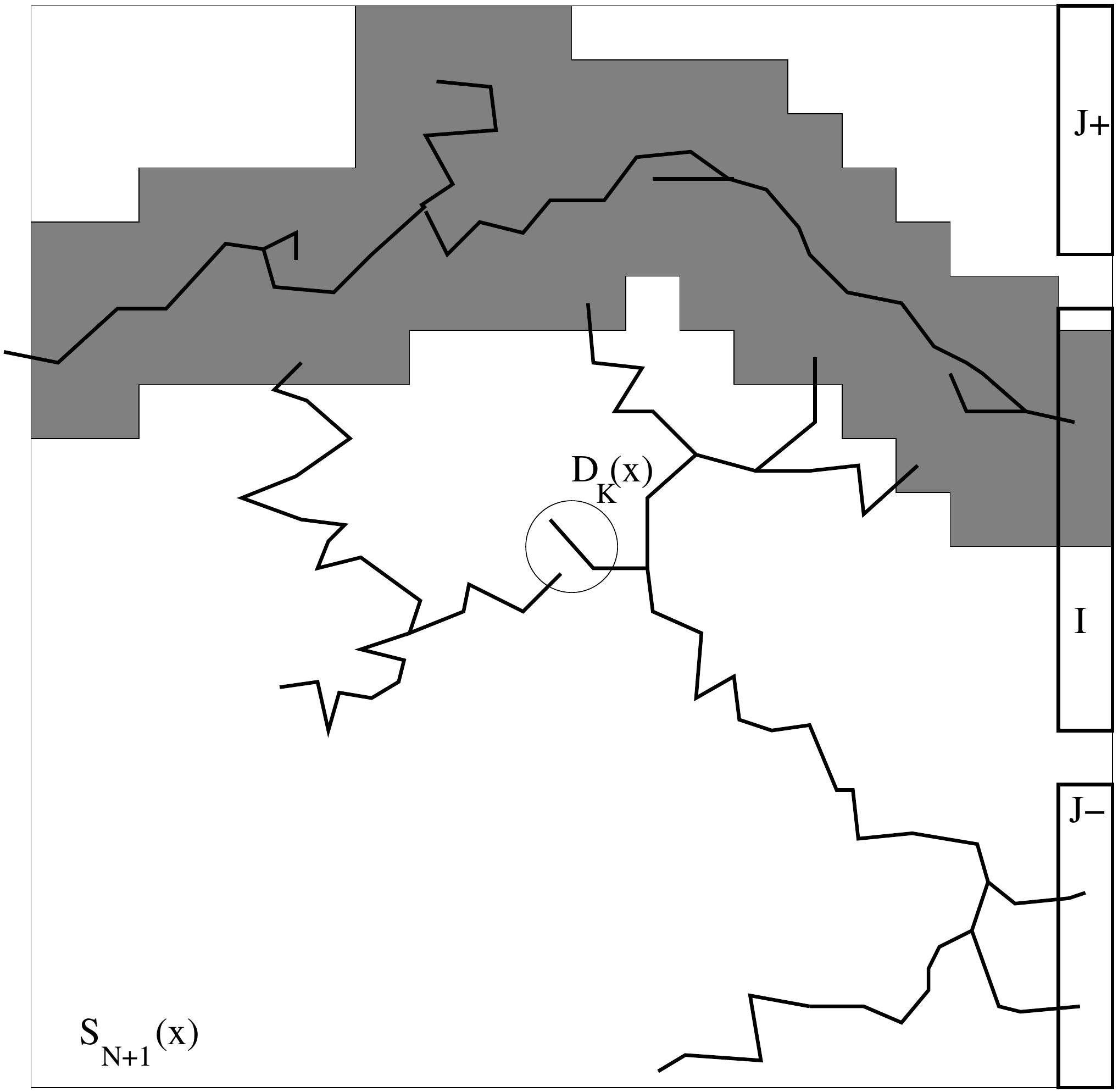}
\caption{The top cluster is the part of $\cC$ within
	  $S_{N+1}(\bx)$. The other clusters are $\cC'$,
	  which is grown within the square $S_{N+1}(\bx)$ but without
	  further exploration from the points created in the shaded region,
	  which is $\cup_{i: \cC \cap Q_i \neq \emptyset} Q_i^+$.
	  In this case $Z=5$; the left-most branch of $\cC'$ reaching
	  the shaded region contributes 2 to $Z$.
	  It is not possible for $\cC'$ to reach both
	  $J^+$ and $J^-$ by this stage.
        }
\end{figure}

  At the next stage
  sample all of the new Poisson points (not part of $\cC$ or $\cC'$) in
  the union of squares $Q_i^{++}$, 
  $i \in {\cal I} $. With probability at least $e^{-25\mu \nu}$
  no new Poisson points are generated at this stage, and if this is the
  case then  $A^+ \cap A^-$ does not occur because the cluster $\cC'$
  dies out without having reached   $J^+$ and $J^-$. 
  Thus we have
  $$
  \Pr[ (A^+ \cap A^-)^c | E \cap \{ Z \leq \nu\} ]
  \geq \exp(-25 \mu \nu)
  \geq \exp(-25 \la \nu)
  .
  $$
  Therefore if  $\Pr[A^+ \cap A^- \cap E] > 1- \delta $ we have
  $$
  \delta > \Pr[ (A^+ \cap A^-)^c \cap E \cap \{ Z \leq \nu\} ]
  > \exp(-25 \la \nu) 
  \Pr[ E \cap \{Z \leq t\}],
  $$
  so that 
$  \Pr[ E \cap \{Z \leq \nu\}] < \eps/3$ by (\ref{deltadef}).
Combined with (\ref{0406a}),
this shows that
\bea
\Pr[(F')^c] \leq \Pr[E^c] + \Pr[(F')^c|E \cap \{Z > \nu\}] + \Pr[E \cap \{Z 
\leq \nu\}]
< \eps.
\eea
 Since
   $(S_{3N+1} \cup S_{N+1}(\bx)) \subset D_{6N} \subset D_{3\|\bx\|} $,
  we have  $F' \subset \{
  D_{K} \leftrightarrow D_K(\bx)$
  in $G(\H_\la \cap D_{3\|\bx\|})\}$.
  By rotation invariance,
  with $M:= \|\bx\|$, 
  we thus have 
  $
  \Pr[F_{K,M,\la} ] >  1- \eps,
  $
  as required.
  \epr

  \subsection{Connection probability}
  \label{subsecCP}
  Given $s>0$, let $V_s,W_s$ be independent uniformly
	  distributed points in
  $B(s)$, independent of $\H_\la$.
  We shall characteraize the giant component of $G(\cH_{\la,s})$
  in terms of those vertices 
  which are path-connected to a fixed disk centred at $\0$. For this, the
  following is useful.
  %In this section we prove the following.
  \begin{Proposition}
	  \label{prop2pt4pt}
	  Suppose 
	  $\la > \la_c(\phi)$. Let $\eps >0$.
	  Then there exists $K>0$ such that
   $\Pr[D_K \leftrightarrow \infty {\rm ~ in~}
   G(\cH_\la) ] > 1-\eps$, and
	  \bea
	  \label{eq2point}
\liminf_{s \to \infty}
	  \Pr[ \{V_s \} \leftrightarrow D_K 
	  ~{\rm in}~ G(\H_{\la,s}^{V_s}
	  ) ] \geq \theta(\phi,\la) - \eps.
	  \eea
	  %Also
	  %\bea
	  %\label{eq4point}
	  %\Pr[ \{V_s \} \leftrightarrow \{W_s\} ~ {\rm and}~
	  %\{X_s\} \leftrightarrow \{Y_s\}
	  %~{\rm in}~ G(\H_{\la,s} \cup\{V_s,W_s,X_s,Y_s\}
	  %) ] \to \theta(\phi,\la)^4.
	  %\eea
  \end{Proposition}
  \bpr
  Assume $\la > \la_c(\phi)$.
  Let $\eps_1 \in (0,\eps/3)$ be chosen such that if $(X_x)_{x \in \Z^2}$ is
  a $7$-dependent Bernoulli random
  field on $\Z^2$ with $\Pr[X_x=1] > 1- 5\eps_1$ for all $x \in \Z^d$,
  than for all $s >0$ and all $x \in \Z^2 \cap B(s)$, there 
  is a lattice path of 1's from $x $ to $\0$ in $\Z^2 \cap B(s)$ with
  probability greater than $1- \eps/3$. 
  The proof that such an $\eps_1$ exists is standard, using
  e.g. \cite{LSS} and a Peierls argument (e.g. \cite[Theorem 9.8]{Penbk}).

  Using Proposition \ref{renormprop}, choose
   $K,M$ such that  $0< K< M/3$, and 
   $$
   \min(\Pr[D_K \leftrightarrow \infty {\rm ~ in~}
   G(\cH_\la) ],
   \Pr[U_{K,M/3,\la}], \Pr[F_{K,M,\la}]) > 1- \eps_1.
   $$

  For each   $ x,y  \in \Z^2$ with $\|x-y\|=1$,
  let 
  $U_x$ denote the event that there is a unique component 
  of $G(\H_{\la} \cap D_{(M/3)+1} ( Mx))$ that meets both $D_K ( Mx)$
  and $\R^2 \setminus  D_{M/3} ( Mx)$.
  Let $F_{xy}:= \{ D_K ( Mx) \leftrightarrow D_K ( My)$ in
  $G(\H_{\la} \cap D_{3M} ( Mx))\}$.

  By translation and rotation invariance of $\H_{\la}$,
  $\Pr[U_x] > 1-\eps_1$ for
  each $x$, and $\Pr[F_{xy}]>1-\eps_1$ for each $(x,y)$.

  For each $x \in \Z^2$,
  let us set $X_x = 1$ if event $U_x$ occurs, and
  also $F_{xy}$ occurs for each of the four $y \in \Z^2$
  with $\|y-x\|=1$;
  otherwise set $X_x =0$. 
   Then 
  % for each $x \in \Z_2$,
   by the union bound $\Pr [X_x=1] \geq 1-5 \eps_1$.
   Also $X_x$ is determined by $\cH_\la|_{D_{3M}(Mx)}$, so
  $(X_x, x \in \Z^2)$ is a $7$-dependent Bernoulli random field.
  By the choice of $\eps_1$,
  for any  $x \in B( s/M ) \cap \Z^2$
  there is a lattice path
  in $\Z^2 \cap  B(s/M ) $ from $x$ to $\0$
  of sites $z $  with $X_z=1$, 
  with probability at least $1-\eps/3$.

  Since $\Pr[ D_{M+K} \leftrightarrow \infty$ in $G(\cH_\la) ]
  \geq
  \Pr[ D_{K} \leftrightarrow \infty$ in $G(\cH_\la)] > 1-\eps_1$, 
   by a similar argument to (\ref{n1eq1}) we can
(and do) take $M_1 >M+K$ such that 
  $\Pr[U_{M+K,M_1,\la} ] > 1- \eps_1$.

  Let $x_V$ be the closest point in $\Z^2 \cap B( s/M )$
  to $M^{-1} V_s$. 
%Then $\|M x_V - V_s\| \leq 2 M$. 
  %and let $Mx_W$ be the closest point in $M (\Z^2 \cap B( s/M ))$
  %to $ W_s$.  
Consider the following events: 
  \begin{itemize}
		  \item
			  $A_1:= \{ 
%is the event that $
V_s \in B(s-4M_1) \setminus D_{3M_1}\}$.
  %and $\|V_s-W_s\| > 3M_1$. 
Provided $s$ is large enough
		  $\Pr[A_1] > 1-\eps/3$.
		  \item
  $ A_2$ is the event that there is a  lattice path from $x_V$
		  to $\0$ within $\Z^2 \cap B(s/M )$
		  with $X_x=1$ for all sites $x$ in the path.
		  By the previous discussion, $\Pr[A_2] > 1- \eps/3$.

	  \item
		  $A_3$ is the event that there is a unique
		  component
		   in $G(\H_{\la}^{V_s} \cap D_{M_1+1} (V_s))$ 
		   that meets both $D_{M+K}(V_s)$ 
		   and $\R^2 \setminus D_{M_1}(V_s)$.
		   Then $\Pr[A_3] > 1- \eps/3$.
		  
	   \item
		   $A_4$ is the event that $\cC_{V_s}(\H_\la^{V_s} ) \cap
		   D_{M_1}(V_s)^c \neq \emptyset$. Then $\Pr[A_4]
		   \geq  \theta(\phi,\la)$.

%	   \item
%		   $A_6$ is the event that $\cC_{W_s}(\H_\la ^{V_s,W_s}) 
%		   \cap
%		   D_{M_1}(V_s)^c \neq \emptyset$.
%		   Then $\Pr[A_6|A_1] \geq \theta(\phi,\la)$.
  \end{itemize}
%Note that
%  %Conditional on $A_1$ occurring, events $A_5$ and $A_6$ are independent,
%  %so
% $\Pr[A_4 |A_1 ] \geq \theta(\phi,\la)$.
%  Therefore $\Pr[A_1 \cap A_4] \geq (1-\eps) \theta(\phi,\la)$.
%  Then 
By the union bound 
  $
  \Pr [ \cap_{i=1}^4 A_i ] \geq \theta(\phi,\la) - \eps,
  $
for all large enough $s$. Therefore it suffices to prove
that if $\cap_{i=1}^4 A_i $
  occurs, then $\{V_s\} \leftrightarrow D_K$ in 
  $G(\H_{\la,s}^{V_s})$.

	  To see this, suppose $\cap_{i=1}^4 A_i$ occurs.
	  Then we have $\|Mx_V-V_s\| \leq M$ so using $A_1$ we have
	  $\dist(Mx_V,B(s)^c) \geq 2 M_1-M > (M/3)+1 $. 
	  Since $X_{x_V} =1$ by $A_2$,
	  there is a unique component of $G(\H_\la \cap D_{(M/3)+1}(Mx_V))$
	  that meets both $D_K(Mx_V)$ and $D_{(M/3)}(Mx_V)^c$, and
	  since $D_{(M/3)+1}(Mx_V) \subset B(s)$, this
	  extends to a unique component of $G(\H_{\la,s})$ that meets
	  both $D_K(Mx_V)$ and $D_{M/3}(Mx_V)^c$. We denote this component by
	  $\cC$.

	  By $A_2$, the component $\cC$ includes a vertex in $D_K$.
	  Choose such a vertex and denote it by $z$.

	  Next, using $A_1$ observe that
	  %$\|V_s\| \geq 3M_1$ by $A_1$, so
  $\|V_s - z\| \geq \|V_s\| - \|z\| \geq 3M_1 -K \geq 2M_1$.
	  Also $D_K(M x_V) \subset D_{M+K}(V_s)$, so $\cC$ meets
	  both  $D_{M+K}(V_s)$ and  $D_{M_1}(V_s)^c$. Therefore
	  using $A_3$ and $A_4$ we have that $V_s $ is connected to $ \cC$,
and thus $\{V_s\} \leftrightarrow D_K$ in $G(\cH_{\la,s}^{V_s})$
as required.
	  \epr

\subsection{Proof of the giant component phenomenon}
\label{sec020111b}

We now write $L_{i,s}$ for $L_i(G(\cH_{\la,s}))$ (we are thinking of
$\la$ and $\phi$ as fixed with $\sup \{r: \phi(r) >0\}=1$). 
For convenience,
we re-state Theorem \ref{thm1}, which we are now ready to prove.
\begin{Theorem}
	\label{thgiantsuper}
	If
	$\la > \la_c(\phi)$, then
	$s^{-2}L_{1,s} \toP \la \theta(\phi,\la)$ and 
	$s^{-2}L_{2,s} \toP 0$ 
	as $s \to \infty$.
\end{Theorem}
\bpr %[Proof of Theorem \ref{thgiantsuper}]
Assume $\la > \la_c(\phi)$.
Let $\eps >0$ and using Proposition \ref{prop2pt4pt}, choose $K>0$ such
that
   $\Pr[D_K \leftrightarrow \infty {\rm ~ in~}
   G(\cH_\la) ] > 1-\eps$, and
(\ref{eq2point}) holds.
Consider the sum
$$
N_s := \sum_{x \in \H_{\la,s}} {\bf 1}\{\{x\} \leftrightarrow D_K
{\rm ~in~} G(\cH_{\la,s}) \}.
$$
%where the sum is over ordered pairs of distinct points of $\cH_{\la,s}$.
Let $V_s,W_s$ be as in Section \ref{subsecCP}. 
%the preceding section. 
By the Mecke formula,
%(\ref{eqMecke1}),
$
\E N_s = \la s^2 
\Pr[\{V_s\} \leftrightarrow D_K {\rm ~ in~} G(\cH_{\la,s} ^{V_s})].
$
%and also 
%$$
%\E[N_s(N_s-1) ] = \la s^2 
%\Pr[\{V_s\} \leftrightarrow D_K
%{\rm ~ in~} G(\cH_{\la,s} \cup\{V_s,W_s,X_s,Y_s\})]
%+ o(s^8).
%$$
Then using (\ref{eq2point}),
%Proposition \ref{prop2pt4pt},
and writing just
$\theta $ for $\theta (\phi, \la)$, we deduce that
\bea
\liminf_{s \to \infty} s^{-2} \E N_s \geq \la (\theta - \eps).
\label{0406b}
\eea

Next, let $N'_s = \sum_{x \in \cH_{\la,s}} {\bf 1}\{ |\cC_x(\cH_{\la,s})|
\geq s^{1/2} \}$ (here $|\cdot|$ represents number of elements).
Using the Mecke formula (\ref{eqMecke1}) we have that
$\E[N'_s]= \la s^2 \Pr[|\cC_{V_s}(\H_{\la,s}^{V_s})|\geq s^{1/2}]$,
and hence
\bea
\lim_{s \to \infty} s^{-2} \E N'_s =  \la \theta. 
\label{0422a}
\eea
Also 
$\E[N'_s(N'_s-1)] = \la^2 s^4
\Pr[ |\cC_{V_s}(\H_{\la,s}^{V_s,W_s})|\geq s^{1/2},
 |\cC_{W_s}(\H_{\la,s}^{V_s,W_s})|\geq s^{1/2}|
]$ 
by (\ref{eqMecke2}),
so that
$$
\lim_{s \to \infty} s^{-4} \E[N'_s(N'_s-1)] = \la^2  \theta^2. 
$$
Thus
$s^{-2} N'_s \to \la \theta$ in $L^2$ and
hence in probability, as $s \to \infty$.

Since $(N_s - N'_s)^+ \leq \H_\la(D_{K+s^{1/2}})$
we have that $s^{-2} \E[(N_s -N'_s)^+] \to 0$ as $s \to \infty$.
Hence by (\ref{0406b}) and (\ref{0422a}),
$$
\limsup \E[s^{-2}(N'_s - N_s)^+] = \limsup \E[s^{-2}(N'_s - N_s)]
\leq \la \eps.
$$
Hence by Markov's inequality
%, as $s \to \infty$ we have
$
\limsup_{s \to \infty}
\Pr[ s^{-2} (N_s'-N_s) \geq \eps^{1/2}] \leq \la \eps^{1/2},
$
and hence
\begin{align}
\limsup_{s \to \infty} \Pr[s^{-2}N_s \leq 
\la  \theta - 2 \eps^{1/2}]
\leq \limsup_{s \to \infty}
	& \left(  \Pr[s^{-2} N'_s \leq \la \theta - \eps^{1/2}]
\right.
\nonumber \\
%~~~~~~~~~~~~~~~~~~~~
	& \left. 
+ \Pr[ s^{-2} (N_s - N'_s) \leq - \eps^{1/2} ] \right)
\leq \la \eps^{1/2}.
\nonumber
\end{align}

As at (\ref{n1eq1}),
we can and do  choose $M_2 > K$ so $\Pr[U_{K,M_2} ] > 1- \eps$.
If $(s/2)> M_2+1$ and  $U_{K,M_2}$ occurs 
then $L_{1,s} \geq N_s - \H_{\la}(D_{M_2})$,
since all $x \in \cH_{\la,s} \setminus D_{M_2}$ that are path-connected
to $D_K$ lie in the same component of $G(\cH_{\la,s})$.
Therefore 
$$
\Pr[s^{-2} L_{1,s}  \leq 
\la  \theta - 3 \eps^{1/2}]
\leq \Pr[s^{-2} N_s \leq \la \theta - 2 \eps^{1/2}]
+\Pr[ s^{-2} \H_{\la}(D_{M_2}) \geq \eps^{1/2}] + \Pr[U_{K,M_2}^c] 
$$
so that
\bea
\limsup_{s \to \infty}
\Pr[s^{-2} L_{1,s}  \leq 
\la  \theta - 3 \eps^{1/2}]
\leq \la \eps^{1/2} + \eps.
\label{0422b}
\eea

Conversely, note that if
$s^2 \la (\theta +\eps) > s^{1/2}$
then
%for $s $ large enough 
$$
\Pr[s^{-2}L_{1,s} \geq \la (\theta + \eps)] \leq
\Pr[ s^{-2} N'_s \geq \la(\theta + \eps) ]
$$  
which tends to zero. Combined with (\ref{0422b})
this shows that
$s^{-2} L_{1,s} \toP \la \theta$.

If 
 $s^2 \la (\theta +\eps) > s^{1/2}$ and 
$L_{1,s} + L_{2,s} \geq s^2 \la (\theta+\eps)$
then either $N'_s \geq s^2 \la (\theta + \eps)$ or
$L_{1,s}+ s^{1/2}\geq s^2 \la (\theta + \eps)$.
Hence
$
\Pr[s^{-2}(L_{1,s}  +L_{2,s}) > \la (\theta + \eps)]
\to 0
.
$
Combined with (\ref{0422b}) this shows that
$s^{-2}( L_{1,s}+ L_{2,s}) \toP \la \theta$
and hence by Slutsky's theorem, 
$s^{-2} L_{2,s} \toP 0$.
\epr


\begin{thebibliography}{}

%	\bibitem{BR}
%		Bollobas and Riordan (2006). {\em Percolation.}
%		Cambridge University Press.
%
%	\bibitem{Grimmett}
%		Grimmett (1999). {\em Percolation. 2nd edition.}
%		Springer.
	%\bibitem{Kingman}
	%	Kingman (1993). {\em Poisson Processes.}

\bibitem{DCST} Duminil-Copin, H., Sidoravicius, V. and Tassion, V. (2016).
	Absence of infinite cluster for critical Bernoulli percolation
		on slabs. {\em Comm. Pure Appl. Math.} {\bf 69},
		1397--1411.

	\bibitem{HLM}
		Heydenreich, M., van der Hofstad, R.,
		Last, G. and  Matzke, K. (2020).
		Lace Expansion and mean-field behavior
for the random connection model. Arxiv:1908.11356.

\bibitem{LP}
	Last, G. and Penrose, M. (2018)
		{\em Lectures on the Poisson process.}
		Cambridge University Press.

	\bibitem{LSS}
		Liggett, T.M., Schonmann, R.H. and Stacey, A.M. (1997).
Domination by product measures. 
		{\em Ann. Probab.} {\bf 25}, 71--95.

	%\bibitem{LPZ}
	%	Last, G. Penrose, M. D. and Zuyev, S. (2017)
	%	On the capacity functional 
	%	of the infinite cluster of a Boolean model.
	%	{\em Ann. Appl. Probab.} {\bf 27} 1678--1701.

	\bibitem{MPS}
		Meester, R., Penrose, M. D. and Sarkar, A. (1997)
The random connection model in high dimensions.
		{\em Statist. Probab. Lett.} {\bf 35}, 145--153. 

	\bibitem{MR} Meester, R. and Roy, R. (1996). 
		{\em Continuum Percolation.}
		Cambridge University Press.

\bibitem{Penbk}
Penrose, M. (2003) {\em Random Geometric Graphs}.
Oxford University Press, Oxford.

\bibitem{P96} Penrose, M.D. (1996)
	Continuum percolation and Euclidean minimal spanning trees in
		high dimensions.
		{\em Ann. Appl. Probab.} {\bf 6}, 528--544. 

\bibitem{P2}
	Penrose, M.D. (2016)  Connectivity of soft random geometric graphs.
		 {\em Ann. Appl. Probab.} {\bf 26},  986--1028. 

	



\end{thebibliography}
\end{document}